\documentclass[a4paper,11pt]{article}

\usepackage[frenchb]{babel}
\usepackage{amsfonts}
\usepackage{amssymb}
\usepackage[T1]{fontenc}

\begin{document}
\newtheorem{prop}{Proposition}
\newtheorem{lemme}{Lemme}
\newtheorem{definition}{D{\'e}finition}
\newtheorem{theorem}{Th{\'e}or{\`e}me}
\newtheorem{coro}{Corollaire}

\author{Rachel Ollivier}

\title{ Modules simples en caract{\'e}ristique $p$ des alg{\`e}bres de Hecke affines de type  $A_2$}
\date{}
\maketitle
{\footnotesize{
Soit $F$ un corps local non archim{\'e}dien de caract{\'e}ristique r{\'e}siduelle $p$ et de corps r{\'e}siduel {\`a} $q$ {\'e}l{\'e}ments.
On ne connait pas  les repr{\'e}sentations irr{\'e}ductibles 
modulo $p$ de $Gl_n(F)$.
Mais il  est tentant de conjecturer que le foncteur des invariants par le pro-$p$-Iwahori
 $I(1)$ les identifie avec les  modules {\`a} droite de l'alg{\`e}bre de Hecke
  ${\cal H}_{\overline{\bf F} _p}(Gl_n(F),I(1))$. 

Pour $n=2$, l'{\'e}tude des ${\cal H}_{\overline{\bf F} _p}(Gl_n(F),I)$-modules simples est le point crucial pour la d{\'e}termination des 
${\cal H}_{\overline{\bf F} _p}(Gl_n(F),I(1))$-modules simples. 
Parmi ceux-l{\`a}, on distingue 
les modules que l'on appelle
"supersinguliers": ils ne correspondent pas, via le foncteur des invariants par le 
pro-$p$-Iwahori, {\`a} des sous-quotients d'induites paraboliques de $Gl_2(F)$.

Nous consid{\'e}rons le cas $n=3$.
Pour d{\'e}terminer les ${\cal H}_{\overline {\mathbf F}_p}(Gl_3(F),I)$-modules simples,
on va  r{\'e}duire les modules standards de 
la ${\overline {\mathbf Q}_p}$-alg{\`e}bre de Hecke-Iwahori  de $Gl_3(F)$ qui admettent
 une structure enti{\`e}re.
On obtient des modules pour la  
${\overline {\mathbf F}_p}$-alg{\`e}bre de Hecke-Iwahori, dont on d{\'e}termine 
les suites de Jordan-H{\"o}lder. On classifie ainsi les modules 
simples de la  ${\overline {\mathbf F}_p}$-alg{\`e}bre de Hecke-Iwahori
 ${\cal H}_{\overline{\bf F} _p}(Gl_3(F),I)$, et l'on d{\'e}crit 
conjecturalement ses modules supersinguliers au sens de [Vig2].}}

\hspace{70mm}

Soit $R$ un anneau commutatif unitaire qui est un ${\bf Z}[q]$-module. On peut consid{\'e}rer $q$ comme une ind{\'e}termin{\'e}e ou bien
l'identifier avec son image $\,q_R$ dans $R$.


\begin{definition}

 Soit $q$ une ind{\'e}termin{\'e}e. On d{\'e}signe par $H$  la $\mathbf{Z}[q]$-alg{\`e}bre engendr{\'e}e par
$T^{\pm 1}, \, S_1$ v{\'e}rifiant les relations
$$(S_1+1)(S_1-q)=0, \:  T^3S_1=S_1T^3,\:
S_1S_2S_1=S_2S_1S_2,\: TS_2=S_1T.$$

\vspace{1mm}

\noindent On d{\'e}signe par  $H_R$ la $R$-alg{\`e}bre $H\otimes _{{\bf Z}[q]} R$.  
Dans $H_R$, on identifie $q$ et son image dans $R$. 
\end{definition}

\vspace{1mm}

\noindent La $R$-alg{\`e}bre de Hecke 
 ${\cal H}_R(Gl_3(F),I)$ est isomorphe {\`a} $H_R$ ([BK] p177).

\vspace{4mm}

\noindent \textbf{Pr{\'e}sentation de Bernstein:} ([Rog1])

\vspace{2mm}

\noindent Soit $A$ l'alg{\`e}bre commutative de la d{\'e}composition de Bernstein
de $H\otimes {\bf Z}[q^{\pm 1}]$,
  isomorphe {\`a} l'alg{\`e}bre des polyn{\^o}mes de Laurent ${\bf Z}[q^{\pm
1}][Y_i^{\pm 1} (1\leq i \leq 3)]$ o{\`u}

\begin{center}$Y_1= q  S_1^{-1} S_2^{-1}T$ ,\,$Y_2=    S_2^{-1} TS_1$,\,$Y_3= q^{-1} T  S_1S_2$.
\end{center}

\noindent L'alg{\`e}bre commutative $A$ est aussi {\'e}gale {\`a}
l'alg{\`e}bre des polyn{\^o}mes de Laurent
\begin{center} ${\bf Z}[q^{\pm 1}][T^{\pm
3}, T_i^{\pm 1} (1\leq i \leq 2)]$ o{\`u} $T_1=T^2S_2S_1, \ \ T_2=TS_1S_2 , $
\end{center} $T_1$, $T_2$ sont les fonctions caract{\'e}ristiques de $Iy_1I$ et de $Iy_2I$ o{\`u}
$y_1=t^2s_2s_1$ est la diagonale $(1,p,p)$ et $y_2=ts_1s_2$ est la
diagonale $(1,1,p)$.

\noindent En effet on a les relations, $Y_1= q  T^{3}  T_1^{-1}, \ \ Y_2=   T_1 T_2^{-1}, \ \ Y_3= q^{-1}
T_2$.

\vspace{2mm}

\noindent D'apr{\`e}s le th{\'e}or{\`e}me de Bernstein,  $A$ et l'alg{\`e}bre de Hecke finie $H_0$ de 
g{\'e}n{\'e}rateurs 
$\,S_1$ et $\,S_2$ engendrent $H\otimes {\bf Z}[q^{\pm 1}]$.

\vspace{4mm}

\noindent Le r{\'e}sultat suivant, relatif {\`a} la structure de l'alg{\`e}bre $H$, est fondamental 
pour la suite. C'est un cas particulier de [Vig2] Th{\'e}or{\`e}me 2. 

\begin{theorem} L'intersection $A\cap H$ est la $\mathbf{Z} [q]$-alg{\`e}bre de type fini
$$\mathbf{Z} [q][\hspace{1mm}(Y_1Y_2Y_3)^{\pm 1},q(Y_iY_j)^{\pm 1}, qY_i^{\pm 1}, 1\leq i\neq j\leq 3\hspace{1mm}].$$
Le centre de $H$ est la $\mathbf{Z} [q]$-alg{\`e}bre de type fini

$$\mathbf{Z} [q][\hspace{1mm}\sigma _3^{\pm 1},\,q \sigma _1, \,q \sigma _2 \hspace{1mm}],$$ o{\`u} $\sigma _1 =Y_1+Y_2+Y_3$, $\sigma _2=Y_1Y_2+Y_2Y_3+Y_1Y_3$
$\sigma _3=Y_1Y_2Y_3$.

\vspace{1mm}

\noindent $H$ est un module de type fini sur son centre.

\end{theorem}

\vspace{3mm}

\noindent Un $R$-caract{\`e}re de $A$ est un morphisme d'anneaux  $\chi: A\rightarrow R$ tel que $\chi (q)=q_R$. Un tel morphisme est enti{\`e}rement d{\'e}termin{\'e} par la donn{\'e}e
du triplet $(y_1,y_2,y_3)$ de ses valeurs en $Y_1,Y_2,Y_3$. 
Ainsi, on a une action naturelle de ${\cal S}_3$ sur les $R$-caract{\`e}res de $A$.

\vspace{4mm}
\noindent \textbf {Rappels sur les modules standards :} ([Rog1])

\vspace{2mm}

Soit $R$ un corps de caract{\'e}ristique nulle. 
Tout  $R$-caract{\`e}re $\chi:A\rightarrow R$ induit un $H_R$-module {\`a} gauche que
l'on appelle $H_R$-module standard associ{\'e} {\`a} $\chi$:
$$I_A(\chi)=H_R\otimes _{A,\chi} R.$$ 
D'apr{\`e}s Bernstein, c'est un  $R$-espace vectoriel de dimension $6$, engendr{\'e} sur $H_R$ 
par le g{\'e}n{\'e}rateur canonique $1\otimes 1$, vecteur propre pour $A$ de valeur propre $\chi$.

\noindent Par adjonction, ces modules standards poss{\`e}dent une propri{\'e}t{\'e} universelle 
fondamentale pour la d{\'e}termination des $H_R$-modules
 simples: 
tout $H_R$-module engendr{\'e} par un vecteur propre pour 
$A$ de valeur propre $\chi$ est quotient du module standard $I_A(\chi)$.

\noindent De plus, deux modules standards $I_A(\chi)$ et $I_A(\chi ')$ 
ont la m{\^e}me suite de Jordan-H{\"o}lder si et seulement si $\chi$ et $\chi '$
 sont 
conjugu{\'e}s sous l'action 
 de ${\cal S}_3$.

\vspace{4mm}

\noindent \textbf{Modules standards entiers:}

\vspace{2mm}

\noindent D{\'e}sormais, $R$ est un anneau commutatif unitaire.

\begin{definition} Soit $\chi:A\cap H\rightarrow R$ un 
$R$-caract{\`e}re de $A\cap H$, c'est {\`a} dire un morphisme d'anneaux tel que $\chi(q)=q_R$.
 Le $H_R$-module standard entier relatif {\`a} $\chi$ est 
le $H_R$-module {\`a} gauche induit de $\chi$
$$I(\chi)=H\otimes_{A\cap H, \chi} R.$$

\end{definition}

\noindent Propri{\'e}t{\'e}s des modules standards entiers:

\noindent(1) Le module standard entier $I(\chi)$ est de type fini sur $R$ car $H$ est de type fini sur $A\cap H$.

\noindent(2) Il est engendr{\'e} par le g{\'e}n{\'e}rateur canonique $\Phi=1\otimes 1$, vecteur propre pour $A\cap H$ de valeur propre $\chi$.
 
\noindent(3) Par adjonction, tout $H_R$-module engendr{\'e}  par un vecteur 
 propre pour $A\cap H$
 de valeur propre $\chi$ est quotient de $I(\chi)$.

\noindent(4) Si $\chi:A\cap H\rightarrow R$ est la restriction de $\chi _A:A\rightarrow R$, alors
$I(\chi)$ est {\'e}gal {\`a} $\,H_R\otimes _{A,\chi _A} R$, le $H_R$-module induit par $\chi _A$.

\vspace{4mm}

\noindent \textbf {Structures enti{\`e}res des modules standards:}

\begin{definition}

Soit $E$ une extension finie de ${\bf Q}_p$ d'anneau d'entiers $O_E$ et 
 $\chi:A\rightarrow E$ un $E$-caract{\`e}re de $A$.

\noindent On dit que $\chi$ est $O_E$-entier si $\chi (A\cap H) \subset O_E$.

\noindent Soit $M$ un   $H_{E}$-module de dimension finie. 
On dit que $M$ a une $O_E$-structure enti{\`e}re
 s' il existe  un sous-$\,O_E$-module de $M$ libre,  stable par $H$ et
qui contient une $E$-base de $M$. On dit alors du $H_E$-module $M$
qu'il est entier.

\noindent On note $\,r_p:O_E\rightarrow \overline{\bf F} _p$ la r{\'e}duction modulo $p$. 
Soit $V$ un $H_{O_E}$-module. On appelle r{\'e}duction de $V$ le $H_{\overline{\bf F}_p}$-module
 $V\otimes _{O_E,r_p}\overline{\bf F}_p$.
\end{definition}
\begin{prop}Soit $E$ une extension finie de ${\bf Q}_p$ d'anneau d'entiers $O_E$ et $\chi_A:A\rightarrow E$ 
un caract{\`e}re de $A$. 

\noindent Le $H_{E}$-module standard $I_A(\chi _A)$  admet 
une $O_E$-structure enti{\`e}re  si et seulement si $\chi _A$ est 
$O_E$-entier,

\noindent c'est {\`a} dire

\begin{center}
$(y_1y_2y_3)^{\pm 1}$, $q(y_iy_j)^{\pm  1}$, $qy_i^{\pm 1}$, $1\leq i\neq j\leq 3$
\end{center}
sont  des {\'e}l{\'e}ments de $O_E$.

\noindent Cette relation est invariante sous l'action de ${\cal S}_3$.

\end{prop}

\noindent \textbf{Preuve:}
C'est [Vig2] (2.5.4).
Soit $\Phi$  
 le g{\'e}n{\'e}rateur canonique de $I_A(\chi _A)$ qui est un vecteur $A$-propre  pour la valeur propre  
$\chi _A$.
 Soit $H_{O_E}\Phi$ le sous $H_{O_E}$-module de $I_A(\chi_A)$ 
engendr{\'e} par $\Phi$.
Plus pr{\'e}cis{\'e}ment, on montre que lorsque $\chi _A$ est  $O_E$-entier sur $A\cap H$, 
le $H_{O_{E}}$-module  $H_{O_E}\Phi$ est une structure enti{\`e}re de 
$I_A(\chi _A)$ appel{\'e}e 
la $O_E$-structure enti{\`e}re canonique.

\vspace{2mm}

\noindent {\textbf{Remarque 2:}} Soit $\chi:A\cap H\rightarrow O_E$ la birestriction 
de $\chi _A$ {\`a} $A\cap H$ et $O_E$.
 Le $\overline {\bf F} _p$-caract{\`e}re 
$r_p\chi:A\cap H\rightarrow \overline{\bf F} _p$ 
induit un $H_{\overline {\bf F} _p}$-module standard entier $I(r_p\chi)$ qui 
est isomorphe {\`a} la r{\'e}duction de la 
structure enti{\`e}re canonique de $I_A(\chi _A)$ ([Vig2] Th{\'e}or{\`e}me 5).

\begin{prop}

Soit $\chi _A: A\rightarrow \overline {\mathbf Q} _p$ un caract{\`e}re. Il existe une extension finie $E$ de ${\mathbf Q}_p$ 
telle que $\chi _A (A)\subset  E $. 

\noindent Si le $H_E$-module standard
$I_A(\chi _A)$ admet une $O_E$-structure enti{\`e}re $L(\chi _A)$, 
la semi-simplification du $H_{\overline {\mathbf F} _p}$-module 
 $L(\chi _A)\otimes \overline {\mathbf F } _p$ ne d{\'e}pend que de $\chi _A$. 
On appelle cette semi-simplification
la r{\'e}duction de $I_A(\chi _A)$.

\noindent De plus, pour tout $w\in {\cal S}_3$, $I_A(w.\chi _A)$ admet 
aussi une structure enti{\`e}re et les r{\'e}ductions de $I_A(\chi _A)$ et $I_A(w.\chi _A)$ 
sont les m{\^e}mes.

\end{prop}

\noindent \textbf{Preuve:} La proposition se d{\'e}duit des propri{\'e}t{\'e}s des modules standards ([Oll] 4.3.1).

\vspace{4mm}

\noindent \textbf {R{\'e}duction des modules standards $I_A(\chi _A)$ pour $\chi _A$ ordinaire:}

\vspace{2mm}

Soit $\chi_{y,A}:A\rightarrow E$  un $E$-caract{\`e}re
donn{\'e} par le triplet $y=(y_1,y_2,y_3)$ de ses valeurs en $Y_1$, 
$Y_2$, $Y_3$. On le suppose entier sur $A\cap H$. Puisque $y_1y_2y_3$ est une unit{\'e}, 
$\chi_{y, A}$ est dans l'orbite sous ${\cal S}_3$ d'un caract{\`e}re tel 
que $y_1$ et $y_3^{-1}$ sont des entiers.
 Pour de tels caract{\`e}res appel{\'e}s \emph{ordinaires}, on d{\'e}finit explicitement
 une structure enti{\`e}re $L(\chi_{y, A})$.

\vspace{1mm}

Soit $\chi_{y,A}: A\rightarrow E$   ordinaire.  On d{\'e}finit des nombres $q(w,y)$, pour tout $w\in {\cal S}_3$:

\begin{quote}

$q(1,y)=1$

\vspace{2mm}

$\begin{array}{ccccc}

q(s_1,y) &= & 1 & si & y_2^{-1} \in O_E \cr
 & =&y_2^{-1}& si & y_2\in O_E \cr
\end{array}$

\vspace{2mm}

 $\begin{array}{ccccc}
q(s_2,y)&=& y_2 & si & y_2^{-1}\in O_E \cr &=&1 & si & y_2\in
O_E\cr
\end{array}$

\vspace{2mm}

$q(s_1s_2,y)=q^{-1}y_1y_2$

\vspace{2mm}

 $  q(s_2s_1,y)=q^{-1}y_1$

\vspace{2mm}

 $\begin{array}{ccccc} q(s_1s_2s_1,y)&=&q^{-1}y_1y_2 & si
& y_2^{-1}\in O_E \cr
  &=&q^{-1}y_1 & si & y_2\in O_E.\cr
\end{array}$

\end{quote}

\noindent Ils sont non-nuls et leurs inverses sont des entiers.

\noindent Soit le syst{\`e}me de Coxeter $({\cal S}_3, \{s_1,s_2\})$, $\ell:{\cal S}_3\rightarrow {\bf N}$ la longueur associ{\'e}e.
On note $T_{s_i}=S_i$, $i=1,2$, 
et $T_{ww'}=T_wT_{w'}$ pour $w, w'$ tels que $\ell(ww')=\ell(w)+\ell(w')$.

\begin{prop}
Soit $\chi_{y, A}:A\rightarrow E$  ordinaire.

\noindent  Soit $L(\chi_{y,A})$ le sous $O_E$-module libre  de rang $6$ de $I_A(\chi_{y, A})$ 
de base  $\{ q(w,y)
 T_w\Phi\} _{w\in {\cal S}_3}$.

\noindent Alors $L(\chi_{y, A})$ est la $O_E$-structure enti{\`e}re canonique de $I_A(\chi_{y, A})$.
\end{prop}

\noindent \textbf{Preuve:}

\noindent Montrons que $L(\chi_{y, A})$ est un $H_{O_E}$-module.
\noindent On va noter $\Phi _w=q(w,y) T_w\Phi$. En particulier,
$\Phi _1=\Phi$. Alors, on a

\begin{center}
$T\Phi=\Phi _{s_2s_1}$, $T^2\Phi=\Phi _{s_1s_2}$,
$T^3\Phi=y_1y_2y_3\Phi$

$T\Phi_{s_1}=\Phi _{s_2}$, $T^2\Phi_{s_1}=\Phi _{s_1s_2s_1}$,
$T^3\Phi_{s_1}=y_1y_2y_3\Phi_{s_1}$
\end{center}

\noindent Donc $T$ stabilise  $L(\chi_{y, A})$ et l'on peut prendre
${\cal B}=( \Phi,\hspace{1mm}
T\Phi,\hspace{1mm}T^2\Phi,\Phi_{s_1},\hspace{1mm}T\Phi_{s_1},\hspace{1mm}T^2\Phi_{s_1})$
comme $O_E$-base de ce module puisque $y_1y_2y_3$ est une unit{\'e}.

\noindent Posons pour simplifier:
$$z=y_1y_2y_3,$$
\vspace{1mm}
$$a=q(s_1,y)^{-1},\hspace{2mm}b=\frac{q(s_2s_1,y)}{q(s_1s_2s_1,y)},\hspace{2mm}c=\frac{q(s_2,y)}{q(s_1s_2,y)},$$
\vspace{1mm}
$$a'=qa^{-1},\hspace{2mm}b'=qb^{-1},\hspace{2mm} c'=qc^{-1}.$$
\noindent On  v{\'e}rifie que ces nombres sont tous des entiers
de $E$ car $\chi_{y,A}$ est ordinaire.

\noindent Les  matrices de $S_1$ et $T$ dans la base ${\cal B}$
sont donn{\'e}es par:

$$
\begin{array}{ccc} S_1 \pmatrix{ 0&0&0&a'&0&0\cr 0&0&0&0&0&b'\cr
0&0&q-1&0&c&0\cr a&0&0&q-1&0&0\cr 0&0&c'&0&0&0\cr
0&b&0&0&0&q-1\cr}& \vspace{2mm} T&  \pmatrix{0&0&z&0&0&0\cr
1&0&0&0&0&0\cr 0& 1&0&0&0&0\cr 0&0&0&0&0&z\cr 0&0&0&1&0&0\cr
0&0&0&0&1&0\cr }
\end{array}$$
Donc $S_1$ et $T$ stabilisent le 
$O_E$-module $L(\chi_{y, A})$ qui est bien stable par $H_{O_E}$. Par cons{\'e}quent, $L(\chi_{y, A})$ contient $H_{O_E}\Phi$.

R{\'e}ciproquement, il est clair que $\Phi$, $T\Phi$, $T^{2}\Phi$ sont des {\'e}l{\'e}ments de $H_{O_E}\Phi$. 
Dans le cas o{\`u} $q(s_1,y)=1$, l'{\'e}l{\'e}ment $\Phi_{s_1}=S_1\Phi$ est dans $H_{O_E}\Phi$.
Dans le cas o{\`u} $q(s_1,y)=y_2^{-1}$, $\Phi_{s_1}=T^{-1}S_2\Phi$ est un {\'e}l{\'e}ment de
 $H_{O_E}\Phi$. D'o{\`u}, $L(\chi_{y, A})=H_{O_E}\Phi$.

\vspace{2mm}

\begin{prop}
Soit $\bar\chi:A\cap H\rightarrow \overline{\bf F}_p$ un caract{\`e}re. 
Il existe un 
caract{\`e}re ordinaire 
$\chi_A:A\rightarrow \overline{\bf Q} _p$ 
tel que la r{\'e}duction du $H_{\overline{\bf Q}_p}$-module standard $I_A(\chi_A)$ 
et le $H_{\overline{\bf F}_p}$-module standard 
entier $I(\bar \chi)$  ont la m{\^e}me semi-simplification.
\end{prop}

\noindent \textbf {Preuve:}
Soit  $\bar\chi:A\cap H\rightarrow \overline{\bf F}_p$. 
On peut montrer {\`a} la main que $\bar \chi$ se rel{\`e}ve en 
un $O_E$-caract{\`e}re $\chi$ de $A\cap H$, o{\`u} $E$ est une 
extension finie de ${\bf Q}_p$. 
Pour une d{\'e}monstration g{\'e}n{\'e}rale, [Vig 2] (3.3).
 Soit $i_p:O_E\rightarrow E$ l'inclusion canonique.
Le $E$-caract{\`e}re $i_p\chi:A\cap H\rightarrow E$ est dans l'orbite sous ${\cal S}_3$ 
d'un $E$-caract{\`e}re ordinaire $\chi_y$. On conclut gr{\^a}ce {\`a} la proposition $2$ et la remarque $2$.


\vspace{6mm}

\noindent Tout caract{\`e}re $\chi _A:A\rightarrow E$ entier sur $A\cap H$ est dans l'orbite sous ${\cal S}_3$ d'un caract{\`e}re ordinaire de la liste suivante.
D'apr{\`e}s la proposition $4$, cette liste donne acc{\`e}s {\`a} la semi-simplification des $H_{\overline{\bf F} _p}$-modules standards entiers.

\noindent La valuation de l'anneau de valuation discr{\`e}te $O_E$ est d{\'e}sign{\'e}e par $val$.


\begin{description}
\item[${L}(\chi_{y, A})\otimes \overline {\mathbf F} _p$ est irr{\'e}ductible] quand:
\begin{itemize}

\item[(1)] $0<val(y_2)<val(q)$ et $val(y_3)=-val(q)$

\item[(2)] $-val(q)<val(y_2)<0$ et $val(y_1)=val(q)$
\vspace{2mm}
\item[(3)] $(\hspace{1mm}val(y_1),\hspace{1mm}val(y_2),\hspace{1mm}val(y_3)\vspace{1mm})=(\hspace{1mm}val(q),\hspace{1mm}0,\hspace{1mm}-val(q)\hspace{1mm})$ et $qy_3\not\equiv y_2$, $qy_2\not\equiv y_1$



\end{itemize}

\item[$L(\chi_{y, A})\otimes \overline {\mathbf F} _p $ est ind{\'e}composable et admet deux sous-quotients irr{\'e}ductibles de dimension $3$] quand:

\begin{itemize}
\item[(4)] $(\hspace{1mm}val(y_1),\hspace{1mm}val(y_2),\hspace{1mm}val(y_3)\hspace{1mm})=(\hspace{1mm}val(q),\hspace{1mm}0,\hspace{1mm}-val(q)\hspace{1mm})$ et $qy_3\not\equiv y_2$, $qy_2\equiv y_1$

\vspace{2mm}
\item[(5)] $(\hspace{1mm}val(y_1),\hspace{1mm}val(y_2),\hspace{1mm}val(y_3)\hspace{1mm})=(\hspace{1mm}val(q),\hspace{1mm}0,\hspace{1mm}-val(q)\hspace{1mm})$ et $qy_3\equiv y_2$, $qy_2\not\equiv y_1$

\vspace{2mm}

\item [(6)]$(\hspace{1mm}val(y_1),\hspace{1mm}val(y_2),\hspace{1mm}val(y_3)\hspace{1mm})=(\hspace{1mm}0,\hspace{1mm}0,\hspace{1mm}0\hspace{1mm})$
\end{itemize}

\item[${L}(\chi_{y, A})\otimes \overline {\mathbf F} _p $ est d{\'e}compos{\'e} en somme directe de  deux sous-quotients irr{\'e}ductibles de
 dimension $3$] quand:

\begin{itemize}
\item[(7)]$0<val(y_2)\leq -val(y_3)<val(q) $
\item[(8)]$-val(q)<- val(y_1) \leq val(y_2)<0$
\end{itemize}

\item[$L(\chi_{y, A})\otimes \overline {\mathbf F} _p$ est ind{\'e}composable de longueur quatre ]quand:

\begin{itemize}
\item[(9)] $(\hspace{1mm}val(y_1),\hspace{1mm}val(y_2),\hspace{1mm}val(y_3)\vspace{1mm})=(\hspace{1mm}val(q),\hspace{1mm}0,\hspace{1mm}-val(q)\hspace{1mm})$ et $qy_3\equiv y_2$, $qy_2\equiv y_1$

\end{itemize}
\end{description}
\

\noindent La liste des sous-quotients irr{\'e}ductibles des ${L}(\chi_{y, A})\otimes \overline {\mathbf F} _p$ peut-{\^e}tre 
d{\'e}crite {\`a} l'aide des $H_{\overline{\bf F}_p}$-modules ci-dessous.
\vspace{2mm}

\noindent Soient $z,y,y'$ des {\'e}l{\'e}ments non nuls de $\overline{\mathbf{F}}_p$ que
 l'on note ici $R$ pour all{\'e}ger les notations. On d{\'e}finit les modules suivants:
\begin{itemize}

\vspace{2mm}
\item[$\bullet$] $L_6(z,y)=Rv \oplus RTv\oplus RT^2 v\oplus Rw \oplus RTw\oplus RT^2 w$

$S_1(v)=0$, $S_1(Tv)=T^2w$, $S_1(T^2v)=-T^2v$, $S_1(w)=-w$, $S_1(Tw)=\frac{y}{z} T^2v$, $S_1(T^2w)=-T^2w$, et $T^3$ agit par le scalaire $z$.
Ce module est irr{\'e}ductible.

\vspace{1mm}
\item[$\bullet$] $\tilde L_6(z,y)=Rv \oplus RTv\oplus RT^2 v\oplus Rw \oplus RTw\oplus RT^2 w$

$S_1(v)=w$, $S_1(Tv)=0$, $S_1(T^2v)=-T^2v$, $S_1(w)=-w$, 
$S_1(Tw)=\frac{y}{z} T^2v$, $S_1(T^2w)=-T^2w$ et $T^3$ agit par le
scalaire $z$. Ce module est irr{\'e}ductible.

\vspace{1mm}
\item[$\bullet$] $K_6(z,y,y')=Rv \oplus RTv\oplus RT^2 v\oplus Rw \oplus RTw\oplus RT^2 w$

$S_1(v)=w$, $S_1(Tv)=\frac{y}{y'}T^2w$, $S_1(T^2v)=-T^2v$, $S_1(w)=-w$, 
$S_1(Tw)=\frac{y'}{z} T^2v$, $S_1(T^2w)=-T^2w$ et $T^3$ agit par le scalaire $z$.

Ce module est irr{\'e}ductible si et seulement si $y^2\neq y'$ et $y'^2\neq zy$.

\item[$\bullet$ ] $M_3(z,y)=Rv \oplus RTv\oplus RT^2 v$, $S_1v=0$, $S_1Tv=y^{-1}T^2v$, $S_1T^2v=-T^2v$, $T^3v=zv$.

\vspace{1mm}

\noindent Alors ce module est irr{\'e}ductible si et seulement si $y^3\neq z$.

\noindent Si  $y^3= z$, alors $M_3(z,y)$   a pour seul sous-module le caract{\`e}re

$$ M_1(0,y):T\rightarrow y, \hspace{1mm} S_1\rightarrow 0,$$ 

et pour seul quotient le module dimension $2$
$$M_2(y)=Rw\oplus RTw \hspace{1mm} S_1w=-w,\hspace{1mm} S_1Tw=0,\hspace{1mm} T^2w=-y^2w-yTw.$$

\item[$\bullet $ ] $\tilde M_3(z,y)=Rv \oplus RTv\oplus RT^2 v$, $S_1v=-v$, $S_1Tv=-y^{-1}T^2v$, $S_1T^2v=-T^2v$, $T^3v=zv$.

\vspace{1mm}

\noindent Alors ce module est irr{\'e}ductible si et seulement si $y^3\neq z$.

\noindent 
Si  $y^3= z$, alors $\tilde M_3(z,y)$  a pour seul sous-module le module de dimension $2$
$$\tilde M_2(y)=Rw\oplus RTw, \hspace{1mm} S_1w=0,\hspace{1mm} S_1Tw=-Tw,\hspace{1mm} T^2w=-y^2w-yTw,$$

\noindent et pour seul quotient le module de dimension $1$
$$M_1(-1,y):T\rightarrow y, \hspace{1mm} S_1\rightarrow -1$$

\item[$\bullet$] $N_3(z,y)=Rv \oplus RTv\oplus RT^2 v$, $S_1v=-v$, $S_1Tv=0$, $S_1T^2v=-yTv-T^2v$, $T^3v=zv$.

\vspace{1mm}

Alors ce module est irr{\'e}ductible si et seulement si $y^3\neq z$.

\noindent Si  $y^3= z$, alors $M_3(z,y)$  a pour seul sous-module le caract{\`e}re  $M_1(-1,y)$
 et pour seul quotient le module de dimension $\,2$  $\,M_2(y)$.

\item[$\bullet  $] $\tilde N_3(z,y)=Rv \oplus RTv\oplus RT^2 v$, $S_1v=0$, $S_1Tv=0$, $S_1T^2v=yTv-T^2v$, $T^3v=zv$.

\vspace{1mm}

Alors ce module est irr{\'e}ductible si et seulement si $y^3\neq z$.

\noindent Si  $y^3= z$, alors $\tilde N_3(z,y)$  a pour seul sous-module le module de
 dimension $\,2$ $\, \tilde M_2(y)$ 
et pour seul quotient le caract{\`e}re  $\,M_1(0,y)$. 
\vspace{1mm}

\item [$\bullet $] $P_3(z)=Rv \oplus RTv\oplus RT^2 v$, $S_1v=0$, $S_1Tv=0$, $S_1T^2v=-T^2v$, $T^3v=zv$.

\noindent  Il est irr{\'e}ductible.
\vspace{1mm}
\item [$\bullet $] $\tilde P_3(z)=Rv \oplus RTv\oplus RT^2 v$, $S_1v=-v$, $S_1Tv=0$, $S_1T^2v=-T^2v$, $T^3v=zv$.

\noindent Il est irr{\'e}ductible.

\end{itemize}
\vspace{2mm}

\vspace{3mm}

\noindent \textbf{Action du centre:}

Soit $z\in \overline{\bf F}_p^*$, $y,y'\in \overline{\bf F}_p$. On d{\'e}signe par $\omega(z,y,y')$ le caract{\`e}re central d{\'e}fini par:

$$\begin{array}{ccc}
 \sigma _3&\mapsto z\cr
q\sigma _1&\mapsto y\cr
q\sigma _2&\mapsto y'\cr
\end{array}$$

\noindent On reprend la num{\'e}rotation pr{\'e}c{\'e}dente pour d{\'e}crire la semi-simplification 
des  
$\,L(\chi_{y, A})\otimes {\overline{\mathbf{F}}_p}$.

\vspace{3mm}

\noindent On pose $z=r_p(y_1y_2y_3)$, $y=r_p(qy_3)$, $y'=r_p(qy_2y_3)$. Alors, les 
$H_{\overline{\bf F}_p}$-modules 
$L(\chi_{y, A})\otimes {\overline{\mathbf{F}}_p}$ de la liste suivante ont pour caract{\`e}re central $\omega (z,y,y')$.

\vspace{1mm}

\noindent \textbf{ Caract{\`e}res centraux r{\'e}guliers: $yy'\neq 0.$}
\begin{itemize}
\item[(3)]
$L(\chi_{y, A})\otimes {\overline{\mathbf{F}}_p}$ est isomorphe au module irr{\'e}ductible  
$K_6(z,y,y')$. On a $y^2\neq y'$ et $y'^2\neq zy$.

\item[(4)]  $L(\chi_{y, A})\otimes {\overline{\mathbf{F}}_p}$ a pour seul sous-module  
$M_3(z,y'y^{-1})$. 
Le quotient obtenu est $\tilde M_3(z,y'y^{-1})$. 
On a $y'^2=zy$, $y^2\neq y'$.


\item[(5)]$L(\chi_{y, A})\otimes {\overline{\mathbf{F}}_p}$ a pour seul sous-module  
$\tilde N_3(z,y'y^{-1})$. 
Le quotient obtenu est $ N_3(z,y'y^{-1})$. On a $y'=y^2$, $y'^2\neq zy$.

\item[(9)]$L(\chi_{y, A})\otimes {\overline{\mathbf{F}}_p}$ a pour sous-quotients irr{\'e}ductibles
 $M_1(0,y)$, $M_1(-1,y)$,
$M_2(y)$, $\tilde M_2(y)$.
On a  $y^2=y'$, $y'^2=zy$.
La d{\'e}composition  de $L(\chi_{y, A})\otimes {\overline{\mathbf{F}}_p}$ est repr{\'e}sent{\'e}e par le diagramme suivant:
\setlength{\unitlength}{1cm}

\begin{picture}(15,3)
\put(0.1,1.3){$L(\chi_{y, A})\otimes \overline{\mathbf{F}}_p$ }
\put (2.5,1.5){\line(1,0){2}}
\put(2.9,1.7){\scriptsize $M_1(-1,y)$}
\put(4.5,1.5){\line(2,1){1.5}}
\put(4.5,1.5){\line(2,-1){1.5}}
\put(6.03,2.25){\line(2,-1){1.5}}
\put(6.03,0.75){\line(2,1){1.5}}
\put(7.54,1.5){\line(1,0){2}}
\put(7.9,1.7){\scriptsize $M_1(0,y)$}
\put(6.65,0.8){\scriptsize $M_2(y)$}
\put(4.55,2.1){\scriptsize $M_2(y)$}
\put(6.65,2.1){\scriptsize $\tilde M_2(y)$}
\put(4.55,.8){\scriptsize $\tilde M_2(y)$}

\put(9.6,1.3){$0$}
 
\end{picture}
\end {itemize}

\noindent \textbf{ Caract{\`e}res centraux singuliers: $yy'=0$ mais $(y,y')\neq (0,0)$.}
\begin{itemize}
\item[(1)]$L(\chi_{y, A})\otimes {\overline{\mathbf{F}}_p}$ est
  isomorphe au module irr{\'e}ductible $L_6(z,y)$. Ici, $y'=0$.
\item[(2)]$L(\chi_{y, A})\otimes {\overline{\mathbf{F}}_p}$ est 
isomorphe au module irr{\'e}ductible $\tilde L_6(z,y')$. Ici $y=0$.
\end{itemize}
\noindent\textbf {Caract{\`e}res centraux supersinguliers: $y=y'=0$.}
 \begin{itemize}

\item[(6)]$L(\chi_{y, A})\otimes {\overline{\mathbf{F}}_p}$ a pour sous-module $\tilde P_3(z)$
 et pour quotient $P_3(z)$.

\item[(7)(8)]$L(\chi_{y, A})\otimes {\overline{\mathbf{F}}_p}$ est 
somme directe de $\tilde P_3(z)$ et $P_3(z)$.
\end{itemize}

\vspace{2mm}



\noindent On a {\'e}tabli une liste de  $H_{\overline{\bf F}_p}$-modules simples:  $M_1(0,y)$, $M(-1,y)$,
 $M_2(y)$, $\tilde M_2(y)$, $ M_3(z,y)$, $\tilde M_3(z,y)$,
 $ N_3(z,y)$, $\tilde N_3(z,y)$,\, $P_3(z)$, $\tilde P_3(z)$
$L_6(z,y)$, $\tilde L_6(z,y)$, $K_6(z,y,y')$ avec $y^2\neq y'$ $y'^2\neq zy$  et $y^3\neq z$. 
On  note que l'action du centre sur ces modules simples d{\'e}termine les valeurs des param{\`e}tres $z$, $y$, $y'$.
Gr{\^a}ce {\`a} cet argument et {\`a}  l'observation de la trace de $\,S_1$ on montre qu'il n'y a pas d'isomorphisme entre
deux modules simples de dimension $\,6$ (resp. $3$) de la liste pr{\'e}c{\'e}dente. On montre {\`a} la main que $M_2(y)$ et $\tilde M_2(y)$ ne sont pas isomorphes.

\vspace{3mm}


On connait la semi-simplification des $H_{\overline{\bf Q}_p}$-modules 
standards ([Rog2]).
 En la comparant, pour les $\chi_{y,A}$ de la liste, avec la semi-simplification 
de $L(\chi_{y,A})\otimes \overline{\bf F}_p$
on remarque que tous les modules simples contenus dans $L(\chi_{y,A})\otimes \overline{\bf F}_p$
 se rel{\`e}vent: ils sont la r{\'e}duction modulo $p$ d'une structure enti{\`e}re d'un $H_{\overline{\bf Q}_p}$-module simple.
\vspace{5mm}

\noindent On a donc obtenu:

\begin{theorem}
\noindent  Tout $H_{\overline{\mathbf F} _p}$-module simple ayant un caract{\`e}re central est de la forme 

\begin{center}
$M_1(0,y)$, $M(-1,y)$,
 $M_2(y)$, $\tilde M_2(y)$, $K_6(z,y,y')$ avec $y^2\neq y'$ et $y'^2\neq zy$

$ M_3(z,y)$, $\tilde M_3(z,y)$,
 $ N_3(z,y)$, $\tilde N_3(z,y)$,\, avec $y^3\neq z$

$L_6(z,y)$, $\tilde L_6(z,y)$, 

 $P_3(z)$, $\tilde P_3(z)$

\end{center}
avec $z,\,y \, y'\,\in \: \overline{\mathbf F} _p ^*$.

\vspace{2mm}

\noindent Un $H_{\overline{\mathbf F} _p}$-module simple  est la r{\'e}duction d'un $H_{\overline{\bf Q}_p}$-module simple entier.

\end{theorem}

\vspace{10mm}


\noindent{\textbf {Bibliographie:}}

\vspace{5mm}

\begin{description}

\item[{[BK]}]Bushnell, C.J.  Kutzko, P.
The admissible dual of $GL(N)$ via compact open subgroups. 
Annals of Mathematics Studies. 129. Princeton University Press, 1993.
\item[{[Oll]}] Ollivier, R. Modules simples en caract{\'e}ristique p 
des alg{\`e}bres de Hecke affines de type $A_1$, $A_2$.
 DEA Institut de Math{\'e}matiques de Jussieu, Juin 2002.

 \item [{[Rog1]}]Rogawski, J.D. 
On modules over the Hecke algebra of a p-adic group. 
Invent. Math 79, 443-465, 1985. 

\item [{[Rog2]}]Rogawski, J.D.
Representations of GL(n) over a p-adic field with an Iwahori-fixed vector
Isr. J. Math 54, 242-256 ,1986.

 \item [{[Vig1]}]Vign{\'e}ras, M.-F   Representations modulo p of the p-adic
 group GL(2,F). Institut de Math{\'e}matiques de Jussieu. Pr{\'e}publication  301, Septembre 2001.

\item[{[Vig2]}] Vign{\'e}ras, M.-F. Alg{\`e}bres de Hecke affines g{\'e}n{\'e}riques. Institut de Math{\'e}matiques de Jussieu. Pr{\'e}publication,  Novembre 2002.

\end{description}

\end{document}